\documentclass[a4paper,10pt]{amsart}%
\usepackage{eurosym}
\usepackage{amsmath}
\usepackage{mathrsfs}
\usepackage{amsfonts}
\usepackage{graphicx}
\usepackage{color}
\usepackage{amsfonts}
\usepackage{amssymb}%

\usepackage{esint}

\usepackage{combelow}
\usepackage{palatino}

\usepackage[abbrev, backrefs]{amsrefs}

\newtheorem{theorem}{Theorem}[section]

\newtheorem{corollary}[theorem]{Corollary}

\newtheorem{definition}[theorem]{Definition}

\newtheorem{lemma}[theorem]{Lemma}

\newtheorem{proposition}[theorem]{Proposition}

\numberwithin{equation}{section}

\newcommand{\R}{\mathbb{R}}

\mathchardef\mhyphen="2D

\begin{document}

\title{Endpoint $L^1$ estimates for Hodge systems}

\author[F. Hernandez]{Felipe Hernandez}
\address[F. Hernandez]{Department of Mathematics,
Building 380, Stanford, California 94305}
\email{felipehb@stanford.edu}

\author[B. Rai\cb{t}\u{a}]{Bogdan Rai\cb{t}\u{a}}
\address[B. Rai\cb{t}\u{a}]{Max-Planck-Instiutut f\"ur Mathematik in den Naturwissenschaften, Inselstrasse 22, 04103 Leipzig, Germany;\newline
Ennio De Giorgi Mathematical Research Center, Scuola Normale Superiore, Piazza dei Cavalieri 7, 56126 Pisa, Italy}
\email{bogdanraita@gmail.com}

\author[D. Spector]{Daniel Spector}
\address[D. Spector]{Department of Mathematics, National Taiwan Normal University, No. 88, Section 4, Tingzhou Road, Wenshan District, Taipei City, Taiwan 116, R.O.C.;\newline
Okinawa Institute of Science and Technology Graduate University,
Nonlinear Analysis Unit, 1919--1 Tancha, Onna-son, Kunigami-gun,
Okinawa, Japan
}
\email{spectda@protonmail.com}

\date{}

\begin{abstract}
	In this paper we give a simple proof of the endpoint Besov-Lorentz estimate
	$$
	\|I_\alpha F\|_{\dot{B}^{0,1}_{d/(d-\alpha),1}(\mathbb{R}^d;\mathbb{R}^k)} \leq C \|F \|_{L^1(\mathbb{R}^d;\mathbb{R}^k)}
	$$
for all $F \in L^1(\mathbb{R}^d;\mathbb{R}^k)$ which satisfy a first order cocancelling differential constraint.  We show how this implies endpoint Besov-Lorentz estimates for Hodge systems with $L^1$ data via fractional integration for exterior derivatives.
\end{abstract}

\maketitle

\section{Introduction}
In the $L^1$ theory for linear elliptic systems it is quite difficult to obtain better than weak-type bounds. A program in this direction was pioneered in the seminal work of J. Bourgain and H. Brezis \cite{BourgainBrezis2007} (see also \cite{BourgainBrezisMironescu,VScurves}) and received remarkable contributions from L. Lanzani and E. Stein \cite{LanzaniStein} and J. Van Schaftingen \cite{VS_divk,VSpams,VS}, while endpoint fine parameter improvements on the Lorentz  \cite{Spector-VanSchaftingen-2018,HS} and Besov-Lorentz \cite{Stolyarov} scales have only recently been obtained.

The purpose of this paper is to give a simple proof of the Besov-Lorentz estimates obtained in \cite{Stolyarov} for a restricted class of operators and to show how this estimate can be used to resolve several open questions in the theory, in particular estimates for Hodge systems \cite[Open Problems 1 \& 2]{VSpams} and the endpoint extension of \cite[Propositions 8.8 \& 8.10]{VS} in the case of first order operators.  Our starting place is an estimate the first and third named authors proved in \cite{HS}, that for any $\alpha \in (0,d)$ there exists a constant $C>0$ for which one has the inequality
\begin{align}
\|I_\alpha F\|_{L^{d/(d-\alpha),1}(\mathbb{R}^d;\mathbb{R}^d)} \leq C \|F \|_{L^1(\mathbb{R}^d;\mathbb{R}^d)}\label{HLS}
\end{align}
for all $F \in L^1(\mathbb{R}^d;\mathbb{R}^d)$ such that $\operatorname*{div}F=0$.  Here $L^{d/(d-\alpha),1}(\mathbb{R}^d;\mathbb{R}^d)$ is a Lorentz space (see Section \ref{Lorentz} for a precise definition) and $I_\alpha$ denotes the Riesz potential of order $\alpha \in (0,d)$, defined for $F \in L^1(\mathbb{R}^d;\mathbb{R}^k)$ by
\begin{align}\label{Riesz_potential}
I_\alpha F(x):= \frac{1}{\Gamma\left(\frac{\alpha}{2}\right)} \int_0^\infty t^{\alpha/2-1} p_t \ast F(x) \;dt,
\end{align}
where $p_t(x):=(4\pi t)^{-d/2} \exp(-|x|^2/4t)$ is the heat kernel in $\mathbb{R}^d$. 

The estimate \eqref{HLS} is a partial replacement for the failure of the Hardy--Littlewood--Sobolev embedding in the $L^1$ endpoint, cf.~\cite[p.~119]{Stein}, while a  comprehensive resolution of the question of a replacement has been given by D. Stolyarov, who in \cite{Stolyarov} (see also \cite[Conjecture 2]{ASW} where such an inequality was conjectured to hold) establishes the sharper inequality
\begin{align}
\|I_\alpha F\|_{\dot{B}^{0,1}_{d/(d-\alpha),1}(\mathbb{R}^d;\mathbb{R}^k)} \leq C \|F \|_{L^1(\mathbb{R}^d;\mathbb{R}^k)}\label{HLS_Besov}
\end{align}
for a very general class of subspaces of $L^1(\mathbb{R}^d;\mathbb{R}^k)$ that includes the kernels of J. Van Schaftingen's class of cocancelling operators \cite{VS} (see Definition \ref{def:cocanceling} below where we recall this class).  The argument in \cite{Stolyarov} is quite involved, and it is there commented by D. Stolyarov that whether the inequality \eqref{HLS_Besov} admits a simpler proof if one only seeks its validity for the more restrictive class of divergence free measures is unknown.  We will shortly give such a proof, which benefited from several insights from his paper and the series of lectures\footnote{We are indebted to D. Stolyarov for the efforts he put into giving these lectures, which can be found at "https://groups.oist.jp/nonlinearanalysis/fall-2020-nonlinear-analysis-seminar-series-0".} he gave on the topic.

To this end, let us recall the approach to \eqref{HLS} in \cite{HS}:  For the space of divergence free measures one finds appropriate atoms, one demonstrates the sufficiency of an estimate on an atom, and one establishes the estimate for a single atom.  The atoms in this case are piecewise-$C^1$ loops which satisfy the uniform ball-growth condition:  To any piecewise-$C^1$ loops $\Gamma \subset \mathbb{R}^d$ one associates the measure which is given by integration along the curve
\begin{align}\label{measure_curve}
\int \Phi\cdot d\mu_\Gamma := \int_0^{|\Gamma|} \Phi(\gamma(t))\cdot \dot{\gamma}(t)\;dt,
\end{align}
for $\Phi \in C_0(\mathbb{R}^d;\mathbb{R}^d)$, where $\gamma :[0,|\Gamma|] \to \mathbb{R}^d$ is the parametrization of $\Gamma$ by arclength.  The atoms are then such piecewise-$C^1$ closed curves for which 
\begin{align}\label{ballgrowth}
\|\mu_\Gamma\|_{\mathcal{M}^1} := \sup_{x \in \mathbb{R}^d, r>0} \frac{||\mu_\Gamma||(B(x,r))}{r} \leq \overline{C}
\end{align}
for some universal constant $\overline{C}>0$, where $||\mu_\Gamma||$ is the total variation measure of $\mu_\Gamma$.  The sufficiency of an estimate on these atoms follows in two steps.  First, by Smirnov's integral decomposition of divergence free measures \cite{Smirnov} one has an approximation of such objects in the strict topology by convex combinations of  $C^1$ closed loops (see \cite{GHS} for the details of this argument).  Second, a surgery on such loops shows how any $C^1$ closed loop $\Gamma$ admits a further decomposition into piecewise-$C^1$ closed loops $\{\Gamma_i\}_{i=1}^{N(\Gamma)}$ which satisfy \eqref{ballgrowth} with some universal constant and whose total length is bounded by a constant times the length of this loop.  This approximation/decomposition and the triangle inequality then yields that it suffices to prove the estimate for a single loop which satisfies the ball growth condition \eqref{ballgrowth}.  Finally, the estimate \eqref{HLS} for a single loop was argued in \cite{HS} by a hands on interpolation that utilizes several pointwise estimates for Riesz potentials and bounds for various maximal functions.

While the argument of \eqref{HLS} in \cite{HS} for a single loop with a ball growth condition involves only estimates for various maximal functions, in this paper we observe that it can be further simplified by the consideration of a very natural \emph{stronger} quantity that arises in Stolyarov's estimates:
\begin{align}\label{Besov_Lorentz_Improvement}
\int_0^\infty t^{\alpha/2-1} \|p_t \ast F\|_{L^{d/(d-\alpha),1}(\mathbb{R}^d;\mathbb{R}^k)} \;dt.
\end{align}
In particular, in \cite{Stolyarov}, Stolyarov shows how if one controls a discrete analogue of \eqref{Besov_Lorentz_Improvement} this implies \eqref{HLS_Besov}.  As we will see below in Section \ref{streamlined}, the continuous version \eqref{Besov_Lorentz_Improvement} also controls the Besov-Lorentz norm and therefore, taking into account the reduction to atoms established in \cite{HS}, for the demonstration of the Besov-Lorentz inequality for divergence free functions it suffices to prove the inequality
\begin{align}\label{HLS_Besov_Lorentz_Improvement}
 \int_0^\infty t^{\alpha/2-1} \|p_t \ast \mu_\Gamma\|_{L^{d/(d-\alpha),1}(\mathbb{R}^d;\mathbb{R}^d)} \;dt \leq C||\mu_\Gamma||(\mathbb{R}^d)
\end{align}
for all piecewise-$C^1$ closed loops $\Gamma$ which satisfy \eqref{ballgrowth}.  Let us remark that it is not difficult to see that \eqref{Besov_Lorentz_Improvement} controls the Lorentz norm of the Riesz potential of a function, since this follows directly from the representation \eqref{Riesz_potential} and Minkowski's inequality for integrals.  The argument for the Besov-Lorentz case is only slightly more complicated because of the more technical definition of the space.

We therefore proceed to argue the validity of the inequality \eqref{HLS_Besov_Lorentz_Improvement}.  We claim this follows easily from the estimates
\begin{align}
 \| p_t\ast \mu_\Gamma\|_{L^{1}(\mathbb{R}^d;\mathbb{R}^d)} &\leq \|p_t\|_{L^1(\mathbb{R}^d)}  ||\mu_\Gamma||(\mathbb{R}^d)
 = ||\mu_\Gamma||(\mathbb{R}^d)
, \label{classical_one}\\
  \| p_t\ast \mu_\Gamma\|_{L^{\infty}(\mathbb{R}^d;\mathbb{R}^d)} &\leq \|p_t\|_{L^\infty(\mathbb{R}^d)} ||\mu_\Gamma||(\mathbb{R}^d)
 = \frac{c}{t^{d/2}}||\mu_\Gamma||(\mathbb{R}^d)
,\label{classical_two}\\
 \| p_t\ast\mu_\Gamma\|_{L^{1}(\mathbb{R}^d;\mathbb{R}^d)} &\leq C_1 \frac{||\mu_\Gamma||(\mathbb{R}^d)^2}{t^{1/2}} \label{one}\\
  \| p_t\ast \mu_\Gamma\|_{L^{\infty}(\mathbb{R}^d;\mathbb{R}^d)} &\leq  \frac{C_2}{t^{(d-1)/2}}\|\mu_\Gamma\|_{\mathcal{M}^1}. \label{two}
\end{align}
The former two are standard (linear) convolution inequalities for $L^1$ functions, while the latter two are \emph{nonlinear} and only hold because we consider closed loops oriented by their tangent.  Indeed, \eqref{one} follows from the fact such objects admit a generalized minimal surface spanning $\Gamma$, while \eqref{two} utilizes the fact that we work with curves (and we later make use of the fact that they satisfy \eqref{ballgrowth}).  Note that if we only utilized \eqref{classical_one} and \eqref{classical_two} it would not be sufficient for our purposes, since for any $1\leq p \leq +\infty$ interpolation would yield the estimate
\begin{align*}
 \| p_t\ast \mu_\Gamma\|_{L^{p,1}(\mathbb{R}^d;\mathbb{R}^d)} \leq  \| p_t\ast \mu_\Gamma\|_{L^{1}(\mathbb{R}^d;\mathbb{R}^d)}^\theta \| p_t\ast \mu_\Gamma\|_{L^{\infty}(\mathbb{R}^d;\mathbb{R}^d)}^{1-\theta},
\end{align*}
where $\theta=1/p$.  In particular, when $p=\frac{d}{d-\alpha}$, using \eqref{classical_one} and \eqref{classical_two} we find
\begin{align*}
 \| p_t\ast \mu_\Gamma\|_{L^{d/(d-\alpha),1}(\mathbb{R}^d;\mathbb{R}^d)} \leq C\frac{||\mu_\Gamma||(\mathbb{R}^d)}{t^{\alpha/2}},
\end{align*}
which is not good enough to get a finite upper bound, as if utilized to estimate the quantity \eqref{Besov_Lorentz_Improvement} gives a logarithmic divergence at both zero and infinity and therefore cannot yield the estimate \eqref{HLS_Besov_Lorentz_Improvement}.

The subtlety is to notice that the combination of \eqref{classical_one} and \eqref{two} gives an estimate with slightly better behavior at zero, while the combination of \eqref{classical_two} and \eqref{one} gives an estimate with slightly better behavior at infinity, the inequalities
\begin{align}
 \| p_t\ast \mu_\Gamma\|_{L^{d/(d-\alpha),1}(\mathbb{R}^d;\mathbb{R}^d)} &\leq C'_1 \frac{||\mu_\Gamma||(\mathbb{R}^d)^{(d-\alpha)/d}}{t^{\alpha(d-1)/2d}},\label{zero} \\
  \| p_t\ast \mu_\Gamma\|_{L^{d/(d-\alpha),1}(\mathbb{R}^d;\mathbb{R}^d)} &\leq C'_2 \frac{||\mu_\Gamma||(\mathbb{R}^d)^{1+(d-\alpha)/d}}{t^{\alpha/2+(d-\alpha)/2d}}.\label{infinity}
\end{align}
Indeed,
\begin{align*}
\alpha(d-1)/2d &<\alpha/2, \\
\alpha/2+(d-\alpha)/2d&>\alpha/2,
\end{align*}
and therefore it remains to divide the integral so as to linearize the estimate, which follows from dividing at $|\Gamma|^2$ (alternatively, one may first reduce to the case $|\Gamma|=1$ by dilation, though we here avoided this argument because the nonlinearity of the estimates \eqref{one} and \eqref{two} becomes less clear).  

We postpone further details until Section \ref{streamlined}, including the proof of the slightly more technical Besov-Lorentz inequality, so that we can continue to a second purpose of this paper,  which is to catalog some implications of the inequality \eqref{HLS_Besov} in the divergence free case.  Indeed, a fundamental contribution of J. Van Schaftingen's paper \cite{VS} is that divergence free vector fields are generic in the class of vector fields which admit a first order cocanceling annihilator. In particular, following his argument we establish
\begin{theorem}\label{mainresult}
Let $L(D)$ be a first order homogeneous linear partial differential operator acting on vector fields $F\colon\R^n\rightarrow\R^k$. Let $\alpha\in(0,d)$. Then the estimate
\begin{align*}
\|I_\alpha F\|_{\dot{B}^{0,1}_{d/(d-\alpha),1}(\mathbb{R}^d;\mathbb{R}^k)} \leq C \|F \|_{L^1(\mathbb{R}^d;\mathbb{R}^k)}\quad\text{for }L(D)F=0
\end{align*}
holds if and only if $L(D)$ is cocanceling, see Definition \ref{def:cocanceling}.
\end{theorem}
We recall that the cocanceling assumption is \emph{very} mild: As was observed in \cite{VS,GR,RaitaSpector}, failure of this assumption is equivalent to the existence of an unconstrained subspace of $L(D)$-free fields.

Beyond an intrinsic interest in  the mapping properties of fractional integrals, the inequality given in Theorem \ref{mainresult} has implications for PDEs.  For example, in \cite{HS} it was demonstrated how \eqref{HLS} implies a Lorentz space sharpening of an estimate of Bourgain and Brezis \cite{BourgainBrezis2004,BourgainBrezis2007}:
If $F \in L^1(\mathbb{R}^3;\mathbb{R}^3)$ is divergence free, the solution of the Div-Curl system
\begin{align*}
\operatorname*{curl} Z &= F   \\
 \operatorname*{div} Z &= 0
\end{align*}
admits the estimate
\begin{align*}
\| Z\|_{L^{3/2,1}(\mathbb{R}^3;\mathbb{R}^3)} \leq  C\| F\|_{L^{1}(\mathbb{R}^3;\mathbb{R}^3)}
\end{align*}
for some $C>0$.  

Theorem \ref{mainresult} of course implies a similar improvement to this inequality, though in this form has easier use in applications.  For example, we immediately obtain
\begin{corollary}\label{mainresult1}
Let $d\geq 2$, $\alpha \in (0,d)$, and $k \in \mathbb{N} \cap [0,d]$.  There exists a constant $C=C(\alpha,d)>0$ such that for $k\leq d-2$
\begin{align}\label{potentialnodiracl1'}
\|I_\alpha d u \|_{\dot{B}^{0,1}_{d/(d-\alpha),1}(\mathbb{R}^d;\Lambda^{k+1}\mathbb{R}^d)} \leq C \| d u\|_{L^1(\mathbb{R}^d;\Lambda^{k+1}\mathbb{R}^d)},
\end{align}
while for $k\geq 2$
\begin{align}\label{potentialnodiracl1''}
\|I_\alpha d^* u \|_{\dot{B}^{0,1}_{d/(d-\alpha),1}(\mathbb{R}^d;\Lambda^{k-1}\mathbb{R}^d)} \leq C \| d^* u\|_{L^1(\mathbb{R}^d;\Lambda^{k-1}\mathbb{R}^d)},
\end{align}
for all $u \in C^\infty_c(\mathbb{R}^d;\Lambda^{k}\mathbb{R}^d)$. 
\end{corollary}
Here, for $k \in \mathbb{N} \cap [0,d]$, $\Lambda^{k}\mathbb{R}^d$ denotes the vector space of $k$-forms, $C^\infty(\mathbb{R}^d;\Lambda^{k}\mathbb{R}^d)$ denotes the space of functions from $\mathbb{R}^d$ to the space of $k$-forms with smooth coefficients, and
\begin{align*}
d &:C^\infty(\mathbb{R}^d;\Lambda^{k}\mathbb{R}^d) \to C^\infty(\mathbb{R}^d;\Lambda^{k+1}\mathbb{R}^d) \\
d^* &:C^\infty(\mathbb{R}^d;\Lambda^{k}\mathbb{R}^d) \to C^\infty(\mathbb{R}^d;\Lambda^{k-1}\mathbb{R}^d)
\end{align*}
are the exterior differential and exterior co-differential, respectively. 

In this form it is easy to obtain improvements to the left-hand-side of estimates for the Hodge system considered by J. Bourgain and H. Brezis in their paper \cite{BourgainBrezis2007} (see also Lanzani and Stein \cite{LanzaniStein} for a simple proof in the spirit of J. Van Schaftingen's simplification \cite{VS} of the original Bourgain-Brezis argument).  In particular, we give an affirmative answer to \cite[Open Problems 1 \& 2]{VSpams}, the following
\begin{theorem}\label{bbq}
Let $d \geq 2$ and $k \in \mathbb{N} \cap [1,d-1]$.  If $F \in L^1(\mathbb{R}^d;\Lambda^{k-1}\mathbb{R}^d)$, and $G \in L^1(\mathbb{R}^d;\Lambda^{k+1}\mathbb{R}^d)$, then the function $Z=d(-\Delta)^{-1} F+d^*(-\Delta)^{-1} G$ satisfies
\begin{align*}
d^* Z &= F   \\
dZ &= G
\end{align*}
and there exists a constant $C>0$ such that
\begin{align*}
\| Z\|_{\dot{B}^{0,1}_{d/(d-1),1}(\mathbb{R}^d;\Lambda^{k}\mathbb{R}^d)} \leq  C\left(\| F\|_{L^1(\mathbb{R}^d;\Lambda^{k-1}\mathbb{R}^d)} +\|G\|_{L^1(\mathbb{R}^d;\Lambda^{k+1}\mathbb{R}^d)}\right),
\end{align*}
where we additionally require $F \equiv 0$ in the case $k=1$ or $G \equiv 0$ in the case $k=d-1$. 
\end{theorem}

Finally, let us record the following duality estimates, which extend \cite[Propositions 8.8 \& 8.10]{VS} to the endpoint $q=\infty$.
\begin{proposition}\label{duality}
	Let $L(D)$ be a first order cocanceling operator on $\R^d$ and $\alpha\in(0,d)$. Then the estimates for vector fields
	\begin{align*}
	\int_{\R^d}F\cdot \varphi\, dx&\leq C\|F\|_{L^1(\R^d;\R^k)}\|D\varphi\|_{L^{d,\infty}(\R^d;\R^{k\times d})},\\
		\int_{\R^d}F\cdot \varphi \,dx&\leq C\|F\|_{L^1(\R^d;\R^k)}\|\varphi\|_{\dot{B}^{\alpha,\infty}_{d/\alpha}(\R^d;\R^k)}
	\end{align*}
	hold if $L(D)F=0$.
\end{proposition}
Here we use a slightly unusual notation for the Besov spaces $\dot{B}^{\alpha,q}_p$, which is consistent with our earlier notation. In other words, $\dot{B}^{\alpha,q}_p=\dot{B}^{\alpha,q}_{p,p}$.

The plan of the paper is as follows.  In Section \ref{streamlined}, we prove Theorem \ref{mainresult}.  In Section \ref{Hodge} we prove Corollary \ref{mainresult1},  Theorem \ref{bbq}, and Proposition \ref{duality}.  In Section \ref{Lorentz} we address an implicit claim in \cite{Spector} that the estimate in the curl free case was optimal on the Lorentz scale.  In particular, we here give a proof of this claim, which in turn, by J. Van Schaftingen's argument implies optimality of the result of the first and third named authors in \cite{HS} on the Lorentz scale.  It is likely these results are optimal on the Besov-Lorentz scale, though we do not have an example which confirms this.  Finally, in Section \ref{curlfree}, we give direct proofs of several of the results for $F \in L^1(\mathbb{R}^d;\mathbb{R}^d)$ such that $\operatorname*{curl}F=0$.  Of course, this is not as general as the divergence free setting, though notably it does not require the surgery construction from \cite{HS} and therefore provides a streamlined proof for the Lorentz inequality that does not require anything beyond the coarea formula and basic interpolation of Lebesgue or Lorentz spaces.

\section{Lorentz and Besov-Lorentz Estimates}\label{streamlined}
We begin with a few more details of the argument of Theorem \ref{mainresult} in the base case $L(D)=\operatorname*{div}$.  Toward the inequality \eqref{HLS}, as discussed in the introduction, the reduction argument given in \cite{HS} implies that it suffices to prove \eqref{Besov_Lorentz_Improvement} for every piecewise $C^1$ closed loop $\Gamma$ that satisfies \eqref{ballgrowth}.  This inequality, in turn, will follow if we can establish the convolution inequalities \eqref{classical_one}, \eqref{classical_two}, \eqref{one}, and \eqref{two}. 

For such curves, \eqref{classical_one} and \eqref{classical_two} follow from standard convolution inequalities, while we now explain in more detail the inequalities \eqref{one} and \eqref{two}.  The inequality \eqref{one} follows from the fact that in Euclidean space piecewise $C^1$ closed loop can be identified with integral currents which admit spanning surfaces.  In particular, by \cite[4.2.10]{Federer}, given $T=\mu_\Gamma$, there exists a (generalized) surface $S$ which satisfies (in a generalized sense)
\begin{align}
\partial S &= T\\
||S||(\mathbb{R}^d) &\leq c ||T||(\mathbb{R}^d)^2. \label{isoperimetric}
\end{align}
From this, one easily argues the estimate \eqref{one} by the computation 
\begin{align*}
 \| p_t\ast \mu_\Gamma\|_{L^{1}(\mathbb{R}^d;\mathbb{R}^d)} &= \|p_t \ast T \|_{L^1(\mathbb{R}^d)} \\
&= \|p_t \ast \partial S \|_{L^1(\mathbb{R}^d)} \\
 &\leq \frac{  \| t^{1/2} |\nabla p_t| \|_{L^1(\mathbb{R}^d)} ||S||(\mathbb{R}^d)}{t^{1/2}},
\end{align*}
the identity
\begin{align*}
 \| t^{1/2} |\nabla p_t| \|_{L^1(\mathbb{R}^d)} = c'
\end{align*}
and the isoperimetric inequality \eqref{isoperimetric}.

Concerning the estimate \eqref{two}, it can be argued even simpler than the $\mathcal{H}^1-BMO$ duality utilized to estimate an analogous quantity in \cite{HS}, as it follows from a simple expansion of the convolution on dyadic annuli, using \eqref{ballgrowth}:  In particular, 
\begin{align*}
|p_t\ast \mu_\Gamma(x)| &\leq  \int_{\mathbb{R}^d} p_t(x-y) \;d||\mu_\Gamma||(y) \\
&= \sum_{n\in \mathbb{Z}} \int_{B(x,2^n\sqrt{t}) \setminus B(x,2^{n-1}\sqrt{t})} p_t(x-y) \;d||\mu_\Gamma||(y) \\
&\leq \sum_{n\in \mathbb{Z}} \frac{1}{(4\pi t)^{d/2}} exp(-2^{2n-2}/4) \int_{B(x,2^n\sqrt{t}) \setminus B(x,2^{n-1}\sqrt{t})}  \;d||\mu_\Gamma||(y) \\
&\leq \sum_{n\in \mathbb{Z}} \frac{1}{(4\pi t)^{d/2}} exp(-2^{2n-2}/4) ||\mu_\Gamma||(B(x,2^n\sqrt{t}))\\
&\leq \sum_{n\in \mathbb{Z}} \frac{C}{(4\pi t)^{d/2}} 2^n\sqrt{t} \;exp(-2^{2n-2}/4) \\
&= \frac{C'}{t^{(d-1)/2}},
\end{align*}
where $C$ is as in \eqref{ballgrowth} and
\begin{align*}
C':=\sum_{n\in \mathbb{Z}} \frac{C}{(4\pi )^{d/2}} 2^n exp(-2^{2n-2}/4).
\end{align*}
This and the argument of the introduction completes the proof of the Lorentz inequality in the case $L(D)= \operatorname*{div}$.  

Concerning the Besov-Lorentz inequality, let us recall one definition of the Besov-Lorentz space $\dot{B}^{0,1}_{d/(d-\alpha),1}(\mathbb{R}^d;\mathbb{R}^k)$, see e.g. \cite[Definition 4.1.2]{AH}:
\begin{align*}
\|F\|_{\dot{B}^{0,1}_{d/(d-\alpha),1}(\mathbb{R}^d;\mathbb{R}^k)}:= \sum_{n \in \mathbb{Z}}\|F\ast (\psi_{2^{n+1}}-\psi_{2^n})\|_{L^{d/(d-\alpha),1}(\mathbb{R}^d;\mathbb{R}^k)},
\end{align*}
where 
\begin{align}\label{dilation}
\psi_r(x)= r^d \psi\left(r x\right)
\end{align}
are dilates of some function $\psi \in \mathcal{S}(\mathbb{R}^d)$ which satisfies 
\begin{align*}
\operatorname*{supp} \widehat{\psi} &\subset B(0,1), \\
\widehat{\psi} \equiv 1 &\text{ on }  B(0,1/2).
\end{align*}

Utilizing this definition, we find that we must estimate
\begin{align*}
\|I_\alpha F\|_{\dot{B}^{0,1}_{d/(d-\alpha),1}(\mathbb{R}^d;\mathbb{R}^k)}= \sum_{n \in \mathbb{Z}}\|I_{\alpha}F\ast (\psi_{2^{n+1}}-\psi_{2^n})\|_{L^{d/(d-\alpha),1}(\mathbb{R}^d;\mathbb{R}^k)}.
\end{align*}
In this form, an observation analogous to that of Stolyarov is that if we define the multiplier
\begin{align*}
\widehat{m}(\xi):=\frac{ \widehat{\psi}(\xi)- \widehat{\psi}(2^{-1}\xi)}{(2\pi|\xi|)^\alpha exp(-4\pi^2 |\xi|^2)},
\end{align*}
then, with the use of the notation for scaling introduced in \eqref{dilation}, one has
\begin{align*}
I_{\alpha}F\ast (\psi_{2^{n+1}}-\psi_{2^n}) = 2^{-n \alpha} p_{2^{-2n}} \ast  F \ast m_{2^{n}}.
\end{align*}
As $\widehat{m}$ is smooth and compactly supported, $m \in L^1(\mathbb{R}^d)$, a space which is invariant with respect to the dilation \eqref{dilation}, i.e.
\begin{align*}
\|m_{2^{n}}\|_{L^1(\mathbb{R}^d)}=\|m\|_{L^1(\mathbb{R}^d)}=c.
\end{align*}
 Therefore, by Young's inequality we have
\begin{align*}
\|I_{\alpha}F\ast (\psi_{2^{n+1}}-\psi_{2^n})\|_{L^{d/(d-\alpha),1}(\mathbb{R}^d;\mathbb{R}^k)} \leq c 2^{-n \alpha} \|p_{2^{-2n}} \ast  F\|_{L^{d/(d-\alpha),1}(\mathbb{R}^d;\mathbb{R}^k)},
\end{align*}
so that summation over $n \in \mathbb{Z}$ gives the inequality
\begin{align}\label{Stolyarov_bound}
\|I_\alpha F\|_{\dot{B}^{0,1}_{d/(d-\alpha),1}(\mathbb{R}^d;\mathbb{R}^k)} &\leq \sum_{n \in \mathbb{Z}} c 2^{-n \alpha} \|p_{2^{-2n}} \ast  F\|_{L^{d/(d-\alpha),1}(\mathbb{R}^d;\mathbb{R}^k)}.
\end{align}

The right-hand-side of this inequality is the discrete quantity that Stolyarov obtains an upper bound for in his paper to prove the Besov-Lorentz bound for the general class of subspaces.  To pass to the continuous version, we use the semi-group property of the heat kernel and the fact that it is a contraction on $L^p$:  For each $n \in \mathbb{Z}$ and all $s \in (2^{-2n-2}, 2^{-2n})$,
\begin{align*}
\|p_{2^{-2n}} \ast  F\|_{L^{d/(d-\alpha),1}(\mathbb{R}^d;\mathbb{R}^k)} &=\|p_{2^{-2n}-s}\ast p_s \ast  F\|_{L^{d/(d-\alpha),1}(\mathbb{R}^d;\mathbb{R}^k)} \\
&\leq \| p_{2^{-2n}-s} \ast |p_s \ast  F|\|_{L^{d/(d-\alpha),1}(\mathbb{R}^d;\mathbb{R}^k)}\\
&\leq \| p_s \ast  F\|_{L^{d/(d-\alpha),1}(\mathbb{R}^d;\mathbb{R}^k)}.
\end{align*}
In particular, integration from $s=2^{-2n-2}$ to $2^{-2n}$ with respect to the measure $ds/s$ gives the inequality
\begin{align*}
\|p_{2^{-2n+2}} \ast  F\|_{L^{d/(d-\alpha),1}(\mathbb{R}^d;\mathbb{R}^k)} &\leq \frac{1}{\ln(2)} \int_{2^{-2n-2}}^{2^{-2n}} \| p_s \ast  F\|_{L^{d/(d-\alpha),1}(\mathbb{R}^d;\mathbb{R}^k)} \frac{ds}{s},
\end{align*}
which in combination with \eqref{Stolyarov_bound} yields
\begin{align*}
\|I_\alpha F\|_{\dot{B}^{0,1}_{d/(d-\alpha),1}(\mathbb{R}^d;\mathbb{R}^k)} &\leq \sum_{n \in \mathbb{Z}} c 2^{-n \alpha}\frac{1}{\ln(2)} \int_{2^{-2n-2}}^{2^{-2n}} \| p_s \ast  F\|_{L^{d/(d-\alpha),1}(\mathbb{R}^d;\mathbb{R}^k)} \frac{ds}{s}.
\end{align*}
By further manipulation we obtain
\begin{align*}
\sum_{n \in \mathbb{Z}} &c 2^{-n \alpha}\frac{1}{\ln(2)} \int_{2^{-2n-2}}^{2^{-2n}} \| p_s \ast  F\|_{L^{d/(d-\alpha),1}(\mathbb{R}^d;\mathbb{R}^k)} \frac{ds}{s} \\
&\leq \sum_{n \in \mathbb{Z}} 2c \frac{1}{\ln(2)} \int_{2^{-2n-2}}^{2^{-2n}} s^{\alpha/2-1} \| p_s \ast  F\|_{L^{d/(d-\alpha),1}(\mathbb{R}^d;\mathbb{R}^k)} ds,
\end{align*}
and thus
\begin{align}\label{continuous_control}
\|I_\alpha F\|_{\dot{B}^{0,1}_{d/(d-\alpha),1}(\mathbb{R}^d;\mathbb{R}^k)} \leq 2c \frac{1}{\ln(2)} \int_{0}^{\infty} s^{\alpha/2-1} \| p_s \ast  F\|_{L^{d/(d-\alpha),1}(\mathbb{R}^d;\mathbb{R}^k)} ds.
\end{align}
The inequality \eqref{continuous_control} is exactly the control of the Besov-Lorentz norm by the quantity \eqref{Besov_Lorentz_Improvement} claimed in the introduction.  In particular, by the argument of the introduction and that preceding in this Section, we have established the estimate claimed in Theorem \ref{mainresult} for $F \in L^1(\mathbb{R}^d;\mathbb{R}^d)$ such that $\operatorname*{div}F=0$.

To conclude the proof of Theorem \ref{mainresult}, we follow the argument of J. Van Schaftingen in \cite{VS} that the general case follows by an algebraic reduction. To do this, we first recall a few facts on differential operators. We will work with first order homogeneous linear differential operators with constant coefficients, which can be written as
 $$
 L(D)F=\sum_{i=1}^d L_i\partial_iF=\sum_{i=1}^d \partial_i(L_iF),
 $$ 
 where $L_i\in\operatorname{Lin}(\R^k;\R^l)\simeq\R^{l\times k}$. Due to the usefulness of Fourier transform for linear equations, it is natural to look at the symbol map
 $$
 L(\xi)=\sum_{i=1}^d \xi_iL_i\in \operatorname{Lin}(\R^k;\R^l)\quad\text{for }\xi\in\R^d.
 $$
 We make the simple observation that we can write
 $$
 L(D)F=\operatorname{div}(TF),\quad\text{where }TF=(L_1F\vert L_2F\vert \ldots\vert L_dF),
 $$
so $T\in\operatorname{Lin}(\R^k,\R^{l\times d})$. The divergence of a matrix field is considered row wise.

We recall the definition of cocancellation:
\begin{definition}\label{def:cocanceling}
	An operator $L(D)$ as above is said to be \emph{cocanceling} if and only if
	$$
	\bigcap_{\xi\in\R^d}\ker L(\xi)=\{0\}.
	$$
 \end{definition}
We will show that cocancellation is equivalent with injectivity of the map $T$ defined above. The following lemma should be compared with \cite[Prop.~2.5]{VS} and the proof of \cite[lem.~3.11]{GR}.
\begin{lemma}
	We have that
	$$
	\bigcap_{\xi\in\R^d}\ker L(\xi)=\ker T.
	$$
\end{lemma}
To prove this, note that a vector $F\in\R^k$ lies in the right hand side if and only if
$$
(TF)\xi=0\text{ for all }\xi\in\R^d\iff TF=0,
$$
which completes the proof.

In particular, $L(D)$ is cocanceling if and only if $T$ is left invertible. If this is the case, we can even write an explicit inverse in terms of the adjoint $T^*$ of $T$,
$$
T^\dagger=(T^*T)^{-1}T^*.
$$
We can thus proceed with the proof of the main result. 

\begin{proof}[Conclusion of the proof of Theorem \ref{mainresult}]
	The necessity of cocancellation follows from by plugging in a Dirac mass in the estimate and noting that $I_\alpha\notin L^{d/(d-\alpha)}$.
	
	Conversely, we already proved the desired estimate for $L(D)=\operatorname{div}$. We note that if $L(D)F=0$, we can write $\operatorname{div}(TF)=0$ and $F=T^\dagger TF$, so that
$$
	I_\alpha F=I_\alpha(T^\dagger TF)=T^\dagger I_\alpha(TF),
	$$
	and using the inequality for divergence free measures and the fact that $T^\dagger, T$ are bounded maps on finite dimensional spaces we obtain
\begin{align*}
	\|I_\alpha F\|_{\dot{B}^{0,1}_{d/(d-\alpha),1}(\mathbb{R}^d;\mathbb{R}^k)}&=\|T^\dagger I_\alpha(TF)\|_{\dot{B}^{0,1}_{d/(d-\alpha),1}(\mathbb{R}^d;\mathbb{R}^k)}\\
	&\leq C\|I_\alpha(TF)\|_{\dot{B}^{0,1}_{d/(d-\alpha),1}(\mathbb{R}^d;\mathbb{R}^{l\times d})}\leq C\|TF\|_{L^1(\R^d;\R^{l\times d})}\\
	&\leq C\|F\|_{L^1(\R^d;\R^{k})},
\end{align*}
	which completes the proof.\hfill
\end{proof}

\section{Hodge Systems and Duality Estimates}\label{Hodge}
We first note that Corollary \ref{mainresult1} follows from Theorem \ref{mainresult} since the  $L^1$ vector fields $du$ and $d^*u$ satisfy the first order conditions $d(du)=0$ and $d^*(d^*u)=0$, which are cocanceling for the claimed ranges of $k$. This algebraic fact is elementary to check, see also \cite[Prop.~3.3]{VS} where the cocancellation of $d$ on $\ell$-forms is proved, $\ell\leq d-1$. The claim for $d^*$ follows by duality.

\begin{proof}[Proof of Theorem \ref{bbq}]
We first infer from Corollary \ref{mainresult1} with $\alpha=1$ that
\begin{align}
\|I_1 F\|_{\dot{B}^{0,1}_{d/(d-1),1}(\mathbb{R}^d;\Lambda^{k-1}\mathbb{R}^d)} \leq C \| F\|_{L^1(\mathbb{R}^d;\Lambda^{k-1}\mathbb{R}^d)}, \label{eq:1}\\
\|I_1 G\|_{\dot{B}^{0,1}_{d/(d-1),1}(\mathbb{R}^d;\Lambda^{k+1}\mathbb{R}^d)} \leq   C\|G\|_{L^1(\mathbb{R}^d;\Lambda^{k+1}\mathbb{R}^d)},\label{eq:1'}
\end{align}
where $F\equiv 0$ if $k=1$ and $G\equiv 0$ if $k=d-1$.  Next, we note that since the Hodge Laplacian coincides with the real variable Laplacian, we can express
\begin{align}\label{eq:2}
 Z=d I_2 F+d^*I_2G=(dI_1)I_1F+(d^*I_1)I_1G,
\end{align}
where we used $(-\Delta)^{-1}=I_2$ and the semigroup property of Riesz potentials. Note that $dI_1, d^*I_1$ give rise to a zero-homogeneous Fourier multiplier, hence can be represented as Calder\'on--Zygmund operators, and are therefore bounded on the Besov-Lorentz spaces (Here one should be careful to note that these operators are mappings from one exterior algebra into another.). It follows from \eqref{eq:1}, \eqref{eq:1'}, and \eqref{eq:2} that (and for the convenience of display we remove the notation in the norm which details the images of each map)
\begin{align*}
	\|Z\|_{\dot{B}^{0,1}_{d/(d-1),1}(\mathbb{R}^d)}&\leq C\left(\|(dI_1)I_1F\|_{\dot{B}^{0,1}_{d/(d-1),1}(\mathbb{R}^d)}+\|(d^*I_1)I_1G\|_{\dot{B}^{0,1}_{d/(d-1),1}(\mathbb{R}^d)}\right)\\
	&\leq C\left(\|I_1F\|_{\dot{B}^{0,1}_{d/(d-1),1}(\mathbb{R}^d)}+\|I_1G\|_{\dot{B}^{0,1}_{d/(d-1),1}(\mathbb{R}^d)}\right)\\
&	\leq C\left(\| F\|_{L^1(\mathbb{R}^d)} +\|G\|_{L^1(\mathbb{R}^d)}\right),
\end{align*}
which completes the proof.
\hfill\end{proof}

\begin{proof}[Proof of Proposition \ref{duality}]
Both inequalities follow from our main result, Theorem \ref{mainresult}.	 To prove the first estimate, we observe that the semi-group property of the Riesz potentials and H\"older's inequality on the Lorentz scale implies
\begin{align*}
\int_{\R^d} F\cdot \varphi\,dx &= \int_{\R^d} I_1F\cdot (-\Delta)^{1/2}\varphi\,dx \\
&\leq C\|I_1F\|_{L^{d/(d-1),1}(\R^d;\R^k)}\| R^* D\varphi \|_{L^{d,\infty}(\R^d;\R^{k})},
\end{align*}
where we denote by $R^*$ the adjoint of the Riesz transforms, $R^*=-\operatorname*{div}I_1$, which satisfies the identity $R^* D\varphi =  (-\Delta)^{1/2}\varphi$.  This inequality, Theorem \ref{mainresult}, and the bound
\begin{align*}
\| R^* D\varphi \|_{L^{d,\infty}(\R^d;\R^{k})} \leq C \|D\varphi \|_{L^{d,\infty}(\R^d;\R^{k\times d})}
\end{align*}
then implies the desired result, the last inequality following from the fact that $R^*$ is bounded on the Lorentz spaces.  

In a similar manner we argue the second inequality of the Proposition.  First, by duality we have
\begin{align*}
\int_{\R^d} F\cdot \varphi\,dx&\leq C\|F\|_{\dot B^{-\alpha,1}_{d/(d-\alpha)}(\R^d;\R^k)}\|\varphi\|_{\dot B^{\alpha,\infty} _{d/\alpha}(\R^d;\R^k)}.
\end{align*}
Next, we observe that the definition of $\dot B^{-\alpha,1}_{d/(d-\alpha)}(\R^d;\R^k)$, in analogy with that of the Besov-Lorentz space utilized in Section \ref{streamlined}, is
\begin{align*}
\|F\|_{\dot B^{-\alpha,1}_{d/(d-\alpha)}(\R^d;\R^k)} = \left(\sum_{n \in \mathbb{Z}} \left(2^{-\alpha n} \|F\ast (\psi_{2^{n+1}}-\psi_{2^n})\|_{L^{d/(d-\alpha)}(\mathbb{R}^d;\mathbb{R}^k)}\right)^{d/(d-\alpha)}\right)^{(d-\alpha)/d}.
\end{align*}
However, the embedding of the sequence space $l^{d/(d-\alpha)} \subset l^1$ yields
\begin{align*}
\|F\|_{\dot B^{-\alpha,1}_{d/(d-\alpha)}(\R^d;\R^k)} \leq \sum_{n \in \mathbb{Z}} 2^{-\alpha n} \|F\ast (\psi_{2^{n+1}}-\psi_{2^n})\|_{L^{d/(d-\alpha)}(\mathbb{R}^d;\mathbb{R}^k)}.
\end{align*}
Finally, we observe that argument of Theorem \ref{mainresult} in Section \ref{streamlined} contains within it the estimate
\begin{align*}
\|F\|_{\dot B^{-\alpha,1}_{d/(d-\alpha)}(\R^d;\R^k)} \leq \|F\|_{L^1(\R^d;\R^k)},
\end{align*}
which completes the proof of our second claim.\hfill
\end{proof}
\section{Optimality on the Lorentz Scale}\label{Lorentz}

In a now classical paper on Sobolev embeddings, A. Alvino \cite{alvino} proved (with sharp constant) that one has
\begin{align*}
\|u\|_{L^{d/(d-1),1}(\mathbb{R}^d)} \leq C \|\nabla u \|_{L^1(\mathbb{R}^d;\mathbb{R}^d)}.
\end{align*}
Such an inequality extends to the case $Du$ is a Radon measure, that is, $Du \in M_b(\mathbb{R}^d;\mathbb{R}^d)$ by approximation.

Here we recall that for $1<p<+\infty$ and $0<q<+\infty$ one defines
\begin{align}\label{Lorentz_def}
\|u\|_{L^{p,q}(\mathbb{R}^n)}^q:= p \int_0^\infty \left(t|\{ |u|>t\}|^{1/p}\right)^q\frac{dt}{t}.
\end{align}

The first result of this Section is a construction which will show the optimality of Alvino's result.

\begin{lemma}\label{optimal_alvino}
For every $q<1$, there exists a sequence $\{u_N\}_{N \in \mathbb{N}} \subset BV(\mathbb{R}^d)$ with 
\begin{align*}
|D u_N|(\mathbb{R}^d)\leq C
\end{align*}
independent of $N \in \mathbb{N}$ and 
\begin{align*}
\lim_{N \to \infty} \|u_N\|_{L^{d/(d-1),q}(\mathbb{R}^d)} =+\infty.
\end{align*}
\end{lemma}
\begin{proof}
Define 
\begin{align*}
u_N = \sum_{i=1}^N h_i \chi_{B(0,r_i)}(x)
\end{align*}
where $\{h_i\}_{i \in \mathbb{N}}$ and $r_i\to 0$ will be chosen later such that
\begin{align*}
|D u_N|(\mathbb{R}^d)= \sum_{i=1}^N h_i \omega_d r_i^{d-1} \leq C
\end{align*}
independent of $N$ and
\begin{align*}
\lim_{N\to \infty} \|u_N\|_{L^{d/(d-1),q}(\mathbb{R}^d)} = +\infty.
\end{align*}

To this end, we define $H_0:=0$,
\begin{align*}
H_i := \sum_{j=1}^i h_j
\end{align*}
and compute
\begin{align*}
\|u_N\|_{L^{d/(d-1),q}(\mathbb{R}^d)}^q &= \frac{d}{d-1} \int_0^{H_N} \left(t|\{ |u|>t\}|^{(d-1)/d}\right)^q\frac{dt}{t} \\
&=  \frac{d}{d-1} \sum_{i=0}^{N-1} \int_{H_{i}}^{H_{i+1}} \left(t \omega_d^{(d-1)/d} r_i^{d-1}\right)^q\frac{dt}{t} \\
&=\frac{d}{d-1}\omega_d^{q(d-1)/d} \sum_{i=0}^{N-1} r_i^{(d-1)q} \int_{H_{i}}^{H_{i+1}} t^{q-1}\;dt.
\end{align*}
Since $q<1$, $t\mapsto t^{q-1}$ is decreasing and therefore
\begin{align*}
\int_{H_{i}}^{H_{i+1}} t^{q-1}\;dt \geq H_{i+1}^{q-1}(H_{i+1}-H_{i}) = H_{i+1}^{q-1} h_{i+1}.
\end{align*}
In particular, with a shift of indices we find
\begin{align*}
\|u_N\|_{L^{d/(d-1),q}(\mathbb{R}^d)}^q \geq \frac{d}{d-1}\omega_d^{q(d-1)/d} \sum_{i=1}^{N} r_i^{(n-1)q}H_i^{q-1} h_i.
\end{align*}
Therefore it remains to choose $h_i,r_i$ such that \begin{align*}
\sum_{i=1}^\infty r_i^{(d-1)q} H_i^{q} \frac{ h_i}{H_i} = +\infty
\end{align*}
and recall we must do so in a way the ensures
\begin{align*}
\sum_{i=1}^\infty h_i \omega_d r_i^{d-1} \leq C.
\end{align*}
Choose $h_i=2^i$, so that $H_i=2^{i+1}$.  Thus we now are left to choose $r_i$ such that
\begin{align*}
\sum_{i=1}^\infty r_i^{(d-1)q} H_i^{q} \frac{ h_i}{H_i} = 2^{q-1}\sum_{i=1}^\infty \left(2^i r_i^{(d-1)}\right)^{q} =+\infty
\end{align*}
and
\begin{align*}
\sum_{i=1}^\infty 2^i r_i^{(d-1)} <+\infty.
\end{align*}
But then the choice $2^i r_i^{(d-1)} = \frac{1}{i^{1/q}}$ is sufficient, as $q<1$.\hfill
\end{proof}

By compactness in $BV(\mathbb{R}^d)$ this implies \begin{corollary}\label{optimal_alvino'}
For every $q<1$ and every $M \in \mathbb{N}$, there exists $u \in BV(\mathbb{R}^d)$ such that
\begin{align*}
\|u \|_{L^{d/(d-1),q}(\mathbb{R}^d;\mathbb{R}^d)} >M |D u|(\mathbb{R}^d).
\end{align*}
\end{corollary}

\begin{proof}
Let $q<1$ and define $h_i=2^i$ and $r_i = \left(\frac{1}{i^{1/q}2^i }\right)^{1/(d-1)}$.  Then by the preceding Lemma, we can find a sequence $\{u_N\}$ which satisfies
\begin{align*}
|D u_N|(\mathbb{R}^d)= \sum_{i=1}^N h_i \omega_d r_i^{d-1} \leq C
\end{align*}
and
\begin{align*}
\|u_N\|_{L^{d/(d-1),q}(\mathbb{R}^d)} \geq \left( \frac{d}{d-1}\omega_d^{q(d-1)/d} 2^{q-1}\sum_{i=1}^N \frac{1
}{i}\right)^{1/q}
\end{align*}
In particular, the choice $u_N$ suffices for $N$ sufficiently large.\hfill
\end{proof}

Observe that Lemma \ref{optimal_alvino} and Corollary \ref{optimal_alvino'} imply the optimality on the Lorentz scale of Alvino's result \cite{alvino}, that the second parameter in the Lorentz scale cannot be taken less than $1$.  In this case it is interesting to note that the construction yields, by the Sobolev inequality and the fact that this sequence has compact support, strong convergence in $L^1(\mathbb{R}^d)$ to the function $u$ given by
\begin{align*}
u = \sum_{i=1}^\infty h_i \chi_{B(0,r_i)}(x).
\end{align*}
In particular, compactness shows $u \in BV(\mathbb{R}^d)$, while following the previous argument shows that with this choice of $h_i,r_i$, one has $u \notin L^{d/(d-1),q}(\mathbb{R}^d)$.

For our purposes here it will be useful to observe another consequence of this construction, the following
\begin{lemma}\label{optimal}
For every $r<1$, there exists a sequence $\{u_N\}_{N \in \mathbb{N}} \subset BV(\mathbb{R}^d)$ with 
\begin{align*}
|D u_N|(\mathbb{R}^d) \leq C
\end{align*}
independent of $N \in \mathbb{N}$ and 
\begin{align*}
\lim_{N \to \infty} \|I_\alpha Du_N \|_{L^{d/(d-\alpha),r}(\mathbb{R}^n;\mathbb{R}^n)} = +\infty.
\end{align*}
\end{lemma}

\begin{proof}
We begin with an the inequality for $u$ in terms of potentials and its gradient
\begin{align*}
|u(x)| \leq c I_{1-\alpha} |I_\alpha Du(x)|.
\end{align*}
We will estimate $u$ by a standard potential estimate for $I_{1-\alpha}$.  To this end, we recall an estimate of Hedberg, see e.g. \cite[Proposition 3.1.2 (a)]{AH}which asserts
\begin{align*}
|I_\beta f| \leq \mathcal{M}(f)^{1-\beta p/d} \|f\|^{\beta p/d}_{L^p(\mathbb{R}^d)}.
\end{align*}
The choice $f=|I_\alpha D u(x)|$, $\beta=1-\alpha$, and $p=d/(d-\alpha)$ yields
\begin{align*}
I_{1-\alpha} |I_\alpha D u(x)| \leq \mathcal{M}(|I_\alpha D u(x)|)^{(n-1)/(n-\alpha)} \||I_\alpha D u|\|^{1-(d-1)/(d-\alpha)}_{L^{d/(d-\alpha)}(\mathbb{R}^d)}.
\end{align*}
By the boundedness of the maximal function (see \cite[Theorem 1.4.19]{Grafakos} for the case $q(d-1)/(d-\alpha)<1$) and properties of the Lorentz spaces, this shows that
\begin{align*}
\|u\|_{L^{d/(d-1),q}(\mathbb{R}^d)} \leq \||I_\alpha D u|\|^{(d-1)/(d-\alpha)}_{L^{d/(d-\alpha),q(d-1)/(d-\alpha)}(\mathbb{R}^n)} \||I_\alpha D u|\|^{1-(d-1)/(d-\alpha)}_{L^{d/(d-\alpha)}(\mathbb{R}^d)}.
\end{align*}
By the embedding proved in \cite{SSVS, Spector}, this implies
\begin{align*}
\|u\|_{L^{d/(d-1),q}(\mathbb{R}^d)} \leq \|| I_\alpha D u|\|^{(d-1)/(d-\alpha)}_{L^{d/(d-\alpha),q(d-1)/(d-\alpha)}(\mathbb{R}^d;\mathbb{R}^d)} |D u|(\mathbb{R}^d)^{1-(d-1)/(d-\alpha)}.
\end{align*}
But then for any $r<1$ we may choose
\begin{align*}
q= r \times \frac{d-\alpha}{d-1} <r<1,
\end{align*}
and the construction from Lemma \ref{optimal_alvino} with this choice of $q$ yields the desired sequence.\hfill
\end{proof}

From this we obtain the optimality of Theorem 1.1 in \cite{Spector}, the following
\begin{corollary}
For every $r<1$ and every $M \in \mathbb{N}$, there exists $u \in BV(\mathbb{R}^d)$ such that
\begin{align*}
\|I_\alpha Du \|_{L^{d/(d-\alpha),r}(\mathbb{R}^d;\mathbb{R}^d)} >M |D u|(\mathbb{R}^d).
\end{align*}
\end{corollary}
\begin{proof}
For any $r<1$, we choose the sequence $\{u_N\}$ from Lemma \ref{optimal_alvino} with the corresponding choice
\begin{align*}
q= r \times \frac{d-\alpha}{d-1} <r<1.
\end{align*}
In particular, for any $M \in \mathbb{N}$ we can find $u \in BV(\mathbb{R}^d)$ such that
\begin{align*}
\|u \|_{L^{d/(d-1),q}(\mathbb{R}^d;\mathbb{R}^d)} >M' |D u|(\mathbb{R}^d).
\end{align*}
Meanwhile, the inequality established in Lemma \ref{optimal} asserts
\begin{align*}
\|u\|_{L^{d/(d-1),q}(\mathbb{R}^d)} \leq \|| I_\alpha D u|\|^{(d-1)/(d-\alpha)}_{L^{d/(d-\alpha),q(d-1)/(d-\alpha)}(\mathbb{R}^d;\mathbb{R}^d)} |D u|(\mathbb{R}^d)^{1-(d-1)/(d-\alpha)},
\end{align*}
and therefore
\begin{align*}
M' |D u|(\mathbb{R}^d) < \|| I_\alpha D u|\|^{(d-1)/(d-\alpha)}_{L^{d/(d-\alpha),r}(\mathbb{R}^d;\mathbb{R}^d)} |D u|(\mathbb{R}^d)^{1-(d-1)/(d-\alpha)}.
\end{align*}
Rearrangement of terms then yields the desired result, with an appropriate choice of $M'$.\hfill
\end{proof}

Here we remark that while compactness in $BV$ again allows one to write down the limit
\begin{align*}
u = \sum_{i=1}^\infty h_i \chi_{B(0,r_i)}(x) \in BV(\mathbb{R}^d)
\end{align*}
the fact that $I_\alpha Du \notin L^{d/(d-\alpha),r}(\mathbb{R}^d;\mathbb{R}^d)$ is not obvious in this case.

\section{Simplifications in the Curl Free case}\label{curlfree}
That one has Lorentz and Besov-Lorentz embeddings for $F \in L^1(\mathbb{R}^d;\mathbb{R}^d)$ such that $\operatorname*{curl}F=0$ is contained in Theorem \ref{mainresult}.  We here will give an even more direct proof, which improves upon those given in \cite{Kolyada,KrantzPelosoSpector,Spector,Spector2} and in this less general setting even simplifies some of the argument from Section \ref{streamlined} above.  In particular, the goal of this Section is to establish  
\begin{theorem}\label{kolyada}
Let $d \geq 2$ and $\alpha \in (0,d)$.  There exists a constant $C=C(\alpha,d)>0$ such that
\begin{align*}
\| F\|_{\dot{B}^{-\alpha,1}_{d/(d-\alpha),1}(\mathbb{R}^d;\mathbb{R}^d)} \leq C \|F \|_{L^1(\mathbb{R}^d;\mathbb{R}^d)}
\end{align*}
for all $F \in L^1(\mathbb{R}^d;\mathbb{R}^d)$ such that $\operatorname*{curl}F=0$.
\end{theorem}
Theorem \ref{kolyada} has an interesting history, as it was pointed out to us by D. Stolyarov that one can deduce it from Theorem 4 of V.I. Kolyada's paper \cite[Theorem 4]{Kolyada}.  The third named author was not aware of this during the writing of \cite{Spector}, which gives a different proof of what is unfortunately a slightly weaker result.  We can here rectify this in giving a simpler proof of the same result.

\begin{proof}[Proof of Theorem \ref{kolyada}]
First we claim that it suffices to prove the estimate
\begin{align*}
\int_0^\infty t^{\alpha/2-1} \| p_t \ast D \chi_E\|_{L^{d/(d-\alpha),1}(\mathbb{R}^d;\mathbb{R}^d)} \;dt \leq C|D\chi_E|(\mathbb{R}^d)
\end{align*}
for all $\chi_E \in BV(\mathbb{R}^d)$, the space of functions of bounded variation.  Indeed, the general case then follows from the coarea formula and the argument in Section \ref{streamlined} as to how the left-hand-side controls the desired Besov-Lorentz norm.  Notice that if one only wants to prove the Lorentz inequality \eqref{HLS}, or the weaker Lebesgue inequality, the essential ingredients at the point are only Minkowski's inequality for integrals and the coarea formula.

Toward establishing the desired inequality for sets of finite perimeter, we recall again the classical convolution inequalities
\begin{align}
\|p_t \ast D\chi_E \|_{L^1(\mathbb{R}^d;\mathbb{R}^d)} &\leq  \| p_t \|_{L^1(\mathbb{R}^d)} |D\chi_E|(\mathbb{R}^d),\label{one_old}\\
\|p_t \ast D\chi_E \|_{L^\infty(\mathbb{R}^d;\mathbb{R}^d)} &\leq \|p_t\|_{L^\infty(\mathbb{R}^d)} |D\chi_E|(\mathbb{R}^d). \label{infinity_old}
\end{align}
Interpolation of these inequalities alone would not suffice, and so we require two additional inequalities which are special to characteristic functions of sets.  It is here that the curl free case is much simpler than the divergence free case, as these two inequalities follow immediately from integration by parts and classical convolution estimates:
\begin{align}
\|p_t \ast D\chi_E \|_{L^1(\mathbb{R}^d;\mathbb{R}^d)} &\leq  \|D p_t\|_{L^1(\mathbb{R}^d;\mathbb{R}^d)} \| \chi_E \|_{L^1(\mathbb{R}^d)},\label{one_new} \\
\|p_t \ast D\chi_E \|_{L^\infty(\mathbb{R}^d;\mathbb{R}^d)} &\leq \|Dp_t\|_{L^1(\mathbb{R}^d;\mathbb{R}^d)} \| \chi_E \|_{L^\infty(\mathbb{R}^d)} \label{infinity_new}.
\end{align}
In particular, one does not need to perform surgery, consider a ball growth condition, maximal function estimates, or generalized minimal surfaces.

From here, one interpolates \eqref{one_old} and \eqref{infinity_new} to obtain
\begin{align}\label{lower}
\|p_t \ast D\chi_E \|_{L^{p,1}(\mathbb{R}^d;\mathbb{R}^d)} \leq \frac{C_3}{t^{1/2p'}} |D\chi_E|(\mathbb{R}^d)^{1/p},
\end{align}
where we have used 
\begin{align*}
\| p_t \|_{L^1(\mathbb{R}^d)}=1, \quad  \| \chi_E \|_{L^\infty(\mathbb{R}^d)}=1, \quad \text{ and } \quad \|Dp_t\|_{L^1(\mathbb{R}^d;\mathbb{R}^d)} =\frac{c}{t^{1/2}}.
\end{align*}

In a similar manner, interpolation of \eqref{one_new} and \eqref{infinity_old} yields
\begin{align}\label{upper}
\|p_t \ast D\chi_E \|_{L^{p,1}(\mathbb{R}^d;\mathbb{R}^d)} &\leq  \frac{C_4}{t^{1/2p+d/2p'}} \| \chi_E \|_{L^1(\mathbb{R}^d)}^{1/p}|D\chi_E|(\mathbb{R}^d)^{1-1/p}
\end{align}
where we have used
\begin{align*}
\|p_t\|_{L^\infty(\mathbb{R}^d)}  = \frac{1}{(4\pi t)^{d/2}} \quad \text{ and } \quad \|Dp_t\|_{L^1(\mathbb{R}^d;\mathbb{R}^d)} =\frac{c}{t^{1/2}}.
\end{align*}
As before, we use the fact that in interpolation of Lebesgue spaces, one can improve the second parameter in the interpolation on the Lorentz scale.  One can do even better here than one, though below one the estimate is no longer linear and so the result for the curl free case does not follow from the inequality for characteristics functions of sets. Notice that if one only wants the Lebesgue scale inequality, the interpolation is an exercise in Real Analysis.

Finally, we can make the estimate by splitting the integral in two pieces
\begin{align*}
\int_0^\infty& t^{\alpha/2-1} \| p_t \ast D \chi_E\|_{L^{d/(d-\alpha),1}(\mathbb{R}^d;\mathbb{R}^d)} \;dt \\
&= \int_0^r t^{\alpha/2-1} \| p_t \ast D \chi_E\|_{L^{d/(d-\alpha),1}(\mathbb{R}^d;\mathbb{R}^d)} \;dt \\
&\;\;+ \int_r^\infty t^{\alpha/2-1} \| p_t \ast D \chi_E\|_{L^{d/(d-\alpha),1}(\mathbb{R}^d;\mathbb{R}^d)} \;dt \\
&=: I+II.
\end{align*}

For $I$, we use the interpolated inequality \eqref{lower} with the choice of $p=\frac{d}{d-\alpha}$ to obtain
\begin{align*}
I &\leq \frac{C_3}{\Gamma\left(\frac{\alpha}{2}\right)} |D\chi_E|(\mathbb{R}^d)^{(d-\alpha)/d} \int_0^r t^{\alpha/2-1-\alpha/2d}\;dt \\
&=  C_3'   |D\chi_E|(\mathbb{R}^d)^{(d-\alpha)/d} r^{\alpha/2-\alpha/2d},
\end{align*}
while for $II$ we use the interpolated inequality \eqref{upper}, with the same choice of $p$ to obtain
\begin{align*}
II &\leq \frac{C_4}{\Gamma\left(\frac{\alpha}{2}\right)} \| \chi_E \|_{L^1(\mathbb{R}^d)}^{(d-\alpha)/d} |D\chi_E|(\mathbb{R}^d)^{\alpha/d} \int_r^\infty t^{\alpha/2-1 -1/2p-d/2p'}\;dt \\
&= C_4' \| \chi_E \|_{M_b(\mathbb{R}^d)}^{(d-\alpha)/d} |D\chi_E|(\mathbb{R}^d)^{\alpha/d} r^{-1/2+\alpha/2d}.
\end{align*}
The desired inequality then follows from optimizing in $r$ and the isoperimetric inequality
\begin{align*}
|E|^{1-1/d} \leq c_d|D\chi_E|(\mathbb{R}^d).
\end{align*}  
This concludes the proof, the Section, and the paper.\hfill
\end{proof}

\section*{Acknowledgements}
The authors would like to thank Dima Stolyarov for his comments on an early draft of the manuscript.  F.H. is supported by the Fannie \& John Hertz Foundation.  Part of this work was undertaken while D.S. was visiting the National Center for Theoretical Sciences in Taiwan.  He would like to thank the NCTS for its support and warm hospitality during the visit.

\bibliographystyle{amsplain}

\begin{thebibliography}{99}
\bib{AH}{book}{
   author={Adams, David R.},
   author={Hedberg, Lars Inge},
   title={Function spaces and potential theory},
   series={Grundlehren der Mathematischen Wissenschaften [Fundamental
   Principles of Mathematical Sciences]},
   volume={314},
   publisher={Springer-Verlag, Berlin},
   date={1996},
   pages={xii+366},
   isbn={3-540-57060-8},
   review={\MR{1411441}},
   doi={10.1007/978-3-662-03282-4},
}
	
	
\bib{ASW}{article}{
   author={Ayoush, Rami},
   author={Stolyarov, Dmitriy},
   author={Wojciechowski, Michal},
   title={Sobolev martingales},
   journal={Rev. Mat. Iberoam.},
   volume={ (to appear)},
   date={},
   number={},
   pages={},
   doi={doi: 10.4171/rmi/1224}
}


	
\bib{alvino}{article}{
   author={Alvino, Angelo},
   title={Sulla diseguaglianza di Sobolev in spazi di Lorentz},
   journal={Boll. Un. Mat. Ital. A (5)},
   volume={14},
   date={1977},
   number={1},
   pages={148--156},
}


\bib{BourgainBrezis2004}{article}{
  author={Bourgain, Jean},
  author={Brezis, Ha{\"{\i}}m},
  title={New estimates for the Laplacian, the div-curl, and related Hodge
  systems},
  journal={C. R. Math. Acad. Sci. Paris},
   volume={338},
  date={2004},
  number={7},
  pages={539--543},
}


\bib{BourgainBrezis2007}{article}{
   author={Bourgain, Jean},
   author={Brezis, Ha{\"{\i}}m},
   title={New estimates for elliptic equations and Hodge type systems},
   journal={J. Eur. Math. Soc. (JEMS)},
   volume={9},
  date={2007},
  number={2},
   pages={277--315},
  }

\bib{BourgainBrezisMironescu}{article}{
   author={Bourgain, Jean},
   author={Brezis, Haim},
   author={Mironescu, Petru},
   title={$H^{1/2}$ maps with values into the circle: minimal
   connections, lifting, and the Ginzburg-Landau equation},
   journal={Publ. Math. Inst. Hautes \'Etudes Sci.},
   number={99},
   date={2004},
   pages={1--115},
}

  

\bib{Federer}{book}{
   author={Federer, Herbert},
   title={Geometric measure theory},
   series={Die Grundlehren der mathematischen Wissenschaften, Band 153},
   publisher={Springer-Verlag New York Inc., New York},
   date={1969},
   pages={xiv+676},
   review={\MR{0257325}},
}




\bib{GHS}{article}{
   author={Goodman, Jesse},
   author={Hernandez, Felipe},
   author={Spector, Daniel},
   title={Two approximation results for divergence free vector fields},
   journal={},
   volume={},
   date={},
   number={},
   pages={},
   issn={},
   review={},
   doi={},
}

\bib{Grafakos}{book}{
   author={Grafakos, Loukas},
   title={Classical Fourier analysis},
   series={Graduate Texts in Mathematics},
   volume={249},
   edition={3},
   publisher={Springer, New York},
   date={2014},
   pages={xviii+638},
   isbn={978-1-4939-1193-6},
   isbn={978-1-4939-1194-3},
   review={\MR{3243734}},
   doi={10.1007/978-1-4939-1194-3},
}

\bib{GR}{article}{
	author={Guerra, Andr\'e},
	author={Rai\cb{t}\u{a}, Bogdan},
	title={Quasiconvexity, null Lagrangians and Hardy space integrability under constant rank constraints},
	journal={preprint, arXiv:1909.03923},
	volume={},
	date={},
	number={},
	pages={},
	issn={},
	review={},
	doi={},
}




\bib{HS}{article}{
   author={Hernandez, Felipe}
   author= {Spector, Daniel},
   title={Fractional Integration and Optimal Estimates for Elliptic Systems},
   journal={preprint, arXiv:2008.05639},
   volume={},
   date={},
   number={},
   pages={},
   issn={},
   review={},
   doi={},
}



\bib{Kolyada}{article}{
   author={Kolyada, V. I.},
   title={On the embedding of Sobolev spaces},
   language={Russian, with Russian summary},
   journal={Mat. Zametki},
   volume={54},
   date={1993},
   number={3},
   pages={48--71, 158},
   issn={0025-567X},
   translation={
      journal={Math. Notes},
      volume={54},
      date={1993},
      number={3-4},
      pages={908--922 (1994)},
      issn={0001-4346},
   },
   review={\MR{1248284}},
   doi={10.1007/BF01209556},
}

\bib{KrantzPelosoSpector}{article}{
   author={Krantz, Steven G.},
   author={Peloso, Marco M.},
   author={Spector, Daniel},
   title={Some remarks on $L^1$ embeddings in the subelliptic setting},
   journal={Nonlinear Anal.},
   volume={202},
   date={2021},
   pages={112149},
   issn={0362-546X},
   review={\MR{4156975}},
   doi={10.1016/j.na.2020.112149},
}

\bib{LanzaniStein}{article}{
   author={Lanzani, Loredana},
   author={Stein, Elias M.},
   title={A note on div curl inequalities},
   journal={Math. Res. Lett.},
   volume={12},
   date={2005},
   number={1},
   pages={57--61},
   issn={1073-2780},
   review={\MR{2122730}},
   doi={10.4310/MRL.2005.v12.n1.a6},
}


\bib{RaitaSpector}{article}{
   author={Rai\cb{t}\u{a}, Bogdan},
   author={Spector, Daniel},
   title={A note on estimates for elliptic systems with $L^1$ data},
   journal={C. R. Math. Acad. Sci. Paris},
   volume={357},
   date={2019},
   number={11-12},
   pages={851--857},
   issn={1631-073X},
   review={\MR{4038260}},
   doi={10.1016/j.crma.2019.11.007},
}

\bib{Riesz}{article}{
   author={Riesz, Fr\'{e}d\'{e}ric},
   title={Sur Une Inegalite Integarale},
   journal={J. London Math. Soc.},
   volume={5},
   date={1930},
   number={3},
   pages={162--168},
   review={\MR{1574064}},
   doi={10.1112/jlms/s1-5.3.162},
}

\bib{SSVS}{article}{
   author={Schikorra, Armin},
   author={Spector, Daniel},
   author={Van Schaftingen, Jean},
   title={An $L^1$-type estimate for Riesz potentials},
   journal={Rev. Mat. Iberoam.},
   volume={33},
   date={2017},
   number={1},
   pages={291--303},
   issn={0213-2230},
   review={\MR{3615452}},
   doi={10.4171/RMI/937},
}
\bib{Smirnov}{article}{
   author={Smirnov, S. K.},
   title={Decomposition of solenoidal vector charges into elementary
   solenoids, and the structure of normal one-dimensional flows},
   language={Russian, with Russian summary},
   journal={Algebra i Analiz},
   volume={5},
   date={1993},
   number={4},
   pages={206--238},
   issn={0234-0852},
   translation={
      journal={St. Petersburg Math. J.},
      volume={5},
      date={1994},
      number={4},
      pages={841--867},
      issn={1061-0022},
   },
   review={\MR{1246427}},
}

\bib{Stein}{book}{
   author={Stein, Elias M.},
   title={Singular integrals and differentiability properties of functions},
   series={Princeton Mathematical Series, No. 30},
   publisher={Princeton University Press, Princeton, N.J.},
   date={1970},
   pages={xiv+290},
   review={\MR{0290095}},
}

\bib{Spector}{article}{
   author={Spector, Daniel},
   title={New directions in harmonic analysis on $L^1$},
   journal={Nonlinear Anal.},
   volume={192},
   date={2020},
   pages={111685, 20},
   issn={0362-546X},
   review={\MR{4034690}},
   doi={10.1016/j.na.2019.111685},
}


\bib{Spector1}{article}{
   author={Spector, Daniel},
   title={An optimal Sobolev embedding for $L^1$},
   journal={J. Funct. Anal.},
   volume={279},
   date={2020},
   number={3},
   pages={108559},
   issn={0022-1236},
   review={\MR{4093790}},
   doi={10.1016/j.jfa.2020.108559},
}


\bib{Spector2}{article}{
   author={Spector, Daniel},
   title={A noninequality for the fractional gradient},
   journal={Port. Math.},
   volume={76},
   date={2019},
   number={2},
   pages={153--168},
   issn={0032-5155},
   review={\MR{4065096}},
   doi={10.4171/pm/2031},
}


\bib{Spector-VanSchaftingen-2018}{article}{
   author={Spector, Daniel},
   author={Van Schaftingen, Jean},
   title={Optimal embeddings into Lorentz spaces for some vector
   differential operators via Gagliardo's lemma},
   journal={Atti Accad. Naz. Lincei Rend. Lincei Mat. Appl.},
   volume={30},
   date={2019},
   number={3},
   pages={413--436},
   issn={1120-6330},
   review={\MR{4002205}},
   doi={10.4171/RLM/854},
}



  \bib{Stolyarov}{article}{
   author={Stolyarov, Dmitriy},
   title={Hardy--Littlewood--Sobolev inequality for $p = 1$},
   journal={arXiv:2010.05297},
   volume={},
   date={},
   number={},
   pages={},
}


\bib{VScurves}{article}{
	author={Van Schaftingen, Jean},
		title={A simple proof of an inequality of Bourgain, Brezis and Mironescu},
	journal={Comptes Rendus Mathematique},
	volume={338},
	number={1},
	pages={23--26},
	year={2004},
	publisher={Elsevier}
}

\bib{VS_divk}{article}{
title={Estimates for L1-vector fields under higher-order differential conditions},
author={Van Schaftingen, Jean},
journal={Journal of the European Mathematical Society},
volume={10},
number={4},
pages={867--882},
year={2008}
}

\bib{VSpams}{article}{
title={Limiting fractional and Lorentz space estimates of differential forms},
author={Van Schaftingen, Jean},
journal={Proceedings of the American Mathematical Society},
volume={138},
number={1},
pages={235--240},
year={2010}
}

\bib{VS}{article}{
   author={Van Schaftingen, Jean},
   title={Limiting Sobolev inequalities for vector fields and canceling
   linear differential operators},
   journal={J. Eur. Math. Soc. (JEMS)},
   volume={15},
   date={2013},
   number={3},
   pages={877--921},
   issn={1435-9855},
   review={\MR{3085095}},
   doi={10.4171/JEMS/380},
}


\end{thebibliography}

\end{document}